\documentclass[12pt]{amsart}
\usepackage{amsfonts}
\usepackage{amsmath,amssymb,amsthm}

\newcommand{\R}{\mathbb R}

\newcommand{\ds}{\displaystyle}
\newtheorem{thm}{Theorem}[section]
\newtheorem{cor}[thm]{Corollary}

\newtheorem{prop}[thm]{Proposition}
\theoremstyle{definition}

\theoremstyle{remark}
\newtheorem{rem}[thm]{Remark}

\begin{document}

\title[CANONICAL WEIERSTRASS REPRESENTATION OF MINIMAL AND MAXIMAL SURFACES]
{Canonical Weierstrass Representation of minimal and maximal surfaces in the
three-dimensional Minkowski space}%
\author{Georgi Ganchev}

\address{Bulgarian Academy of Sciences, Institute of Mathematics and Informatics,
Acad. G. Bonchev Str. bl. 8, 1113 Sofia, Bulgaria}%
\email{ganchev@math.bas.bg}%

\subjclass[2000]{Primary 53A35, Secondary 53B30}%
\keywords{Strongly regular time-like surfaces, strongly regular space-like surfaces,
canonical representation of minimal surfaces, canonical representation of maximal surfaces}%

\begin{abstract}

We prove that any minimal (maximal) strongly regular surface in the three-dimensional
Minkowski space locally admits canonical principal parameters. Using this result, we find
a canonical representation of minimal strongly regular time-like surfaces, which makes
more precise the Weierstrass representation and shows more precisely the correspondence
between these surfaces and holomorphic functions (in the Gauss plane). We also find a
canonical representation of maximal strongly regular space-like surfaces, which makes more
precise the Weierstrass representation and shows more precisely the correspondence between
these surfaces and holomorphic functions (in the Lorentz plane). This allows us to describe
locally the solutions of the corresponding natural partial differential equations.
\end{abstract}

\maketitle

\section{Introduction}
In \cite{GM} we proved that any minimal strongly regular surface in Euclidean space can
be endowed locally with canonical principal parameters. Using this fact, in \cite{G} we
found a canonical Weierstrass representation of minimal strongly regular surfaces.

In this paper we consider strongly regular surfaces in the three-dimensional Minkowski
space ${\R}^3_1$.

We prove that any minimal strongly regular time-like surface can be endowed locally
with canonical principal parameters. Using this result we prove the following
\vskip 2mm
\noindent
{\bf Theorem 1. (Canonical Weierstrass representation of minimal time-like surfaces)}
{\it Any minimal strongly regular  time-like surface
$\mathcal M: \, \textbf{z}=\textbf{z}(x,y), \; (x,y) \in {\mathcal D} \subset \mathbb{C}$
parameterized with canonical principal parameters has locally a representation of the type
$$\begin{array}{l}
\ds{(z_1)_x-i(z_1)_y=\;\;\,\frac{1}{2} \, \frac{w^2+1}{w'}},\\
[4mm]
\ds{(z_2)_x-i(z_2)_y=-\frac{i}{2} \, \frac{w^2-1}{w'}},\\
[4mm]
\ds{(z_3)_x-i(z_3)_y=-\frac{w}{w'}},
\end{array}$$
where $w=u(x,y)+iv(x,y)$ is a holomorphic function in $\mathbb{C}$ satisfying the conditions
$$u^2+v^2-1 \neq 0, \qquad \mu:=\frac{(u_x^2+u_y^2)}{(u^2+v^2-1)^2};$$
$$\mu > 0, \quad \mu_x \mu_y \neq 0. $$}

As a consequence of this theorem we obtain a local description of the solutions of the natural partial
differential equation of minimal time-like surfaces.

Further we apply the above scheme for maximal strongly regular space-like surfaces in
${\R}^3_1$.

We prove that any maximal strongly regular space-like surface admits locally canonical principal
parameters. Using this theorem we prove
\vskip 2mm
\noindent
{\bf Theorem 2. (Canonical Weierstrass representation of maximal space-like surfaces)}
{\it Any maximal strongly regular  space-like surface
$\mathcal M: \, \textbf{z}=\textbf{z}(x,y), \; (x,y) \in {\mathcal D} \subset \mathbb{L}$
parameterized with canonical principal parameters has locally a representation of the type
$$\begin{array}{l}
\ds{(z_1)_x+j(z_1)_y=-\frac{w}{w'}}\\
[3mm]
\ds{(z_2)_x+j(z_2)_y=\frac{1}{2} \, \frac{w^2-1}{w'}},\\
[3mm]
\ds{(z_3)_x+j(z_3)_y=\frac{j}{2} \, \frac{w^2+1}{w'}},
\end{array}$$
where $w=u(x,y)+jv(x,y)$ is a holomorphic function in $\mathbb{L}$ satisfying the conditions
$$u^2-v^2+1 \neq 0, \qquad \mu:=\frac{(u_x^2-u_y^2)}{(u^2-v^2+1)^2};$$
$$\mu > 0, \quad \mu_x \mu_y \neq 0.$$}

As a consequence of this theorem we obtain a local description of the solutions of the natural
partial differential equation of maximal space-like surfaces.

\section{Minimal time-like surfaces and canonical principal parameters}

Let $\mathcal M: \, \textbf{z}=\textbf{z}(x,y), \; (x,y) \in {\mathcal D}$ be a space-like
surface in hyperbolic space ${\R}_1^3$ and $\nabla$ be the flat Levi-Civita connection of the
standard metric in ${\R}_1^3$. The unit normal vector field to ${\mathcal M}$ is denoted
by $l$ and $E, F, G; \; e, f, g$ stand for the coefficients of the first and the second
fundamental forms, respectively. In this case we have
$$E=\textbf{z}_x^2>0, \quad G=\textbf{z}_y^2>0, \quad EG-F^2>0, \quad l^2=-1.$$

We suppose that the surface has no umbilical points and the principal lines on
$\mathcal M$ form a parametric net, i.e.
$$F(x,y)=f(x,y)=0, \quad (x,y) \in \mathcal D.$$
Then the principal curvatures $\nu_1, \nu_2$ and the principal geodesic curvatures
(geodesic curvatures of the principal lines) $\gamma_1, \gamma_2$ are given by
$$\nu_1=\frac{e}{E}, \quad \nu_2=\frac{g}{G}; \qquad
\gamma_1=-\frac{E_y}{2E\sqrt G}, \quad \gamma_2= \frac{G_x}{2G\sqrt E}. \leqno(2.1)$$

We consider the tangential frame field $\{X, Y\}$ determined by
$$X:=\frac{\textbf{z}_x}{\sqrt E}, \qquad Y:=\frac{\textbf{z}_y}{\sqrt G}$$
and suppose that the moving frame $XYl$ is positive oriented so that
$$\nu_1-\nu_2>0.$$
Then the following Frenet type formulas for the frame field $XYl$ are valid

$$\begin{tabular}{ll}
$\begin{array}{llccc}
\nabla_{X} \,X & = &  &\gamma_1 \,Y - \nu_1 \, l,  &\\
[2mm]
\nabla_{X} Y & = -\gamma_1 \, X, & & \\
[2mm]
\nabla_{X} \, l & = - \nu_1 \, X, & & &
\end{array}$ &
\quad
$\begin{array}{llccc}
\nabla_{Y} \,X & = & & \gamma_2 \, Y, &\\
[2mm]
\nabla_{Y} Y & = -\gamma_2 \, X & &  & - \nu_2 \, l,\\
[2mm]
\nabla_{Y}\, l & = &  & - \nu_2 \, Y. &
\end{array}$
\end{tabular}\leqno(2.2)$$
\vskip 2mm

The Codazzi equations have the following form
$$\gamma_1=\frac{Y(\nu_1)}{\nu_1-\nu_2}= \frac{(\nu_1)_y}{\sqrt G(\nu_1-\nu_2)},
\qquad
\gamma_2=\frac{X(\nu_2)}{\nu_1-\nu_2}=
\frac{(\nu_2)_x}{\sqrt E(\nu_1-\nu_2)} \leqno(2.3)$$
and the Gauss equation can be written as follows:
$$X(\gamma_2)-Y(\gamma_1)+ \gamma_1^2+\gamma_2^2 = \nu_1\nu_2 = K,$$
or
$$\frac{(\gamma_2)_x}{\sqrt E}-\frac{(\gamma_1)_y}{\sqrt G}+
\gamma_1^2+\gamma_2^2=\nu_1\nu_2=K.$$

A time-like surface ${\mathcal M}: \; \textbf{z}=\textbf{z}(x,y), \; (x,y)\in \mathcal D$
parameterized with principal parameters is said to be \emph{strongly regular} if (cf \cite{GM})
$$\gamma_1(x,y)\gamma_2(x,y) \neq 0, \quad (x,y)\in\mathcal D.$$

The Codazzi equations (2.3) imply that
$$\gamma_1 \gamma_2 \neq 0 \; \iff \; (\nu_1)_y (\nu_2)_x \neq 0.$$
Then the following formulas
$$\sqrt E=\frac{(\nu_2)_x}{\gamma_2(\nu_1-\nu_2)} >0, \quad
\sqrt G=\frac{(\nu_1)_y}{\gamma_1 (\nu_1-\nu_2)}>0 \, \leqno(2.4)$$
are valid for strongly regular surfaces. Because of (2.4) formulas (2.2) become

$$\begin{array}{lccc}
X_x  = &  & \displaystyle{\frac{\gamma_1 \, (\nu_2)_x}
{\gamma_2(\nu_1-\nu_2)}\, Y}
&- \, \displaystyle{\frac{\nu_1 \, (\nu_2)_x}{\gamma_2
(\nu_1-\nu_2)}\, l},\\
[4mm]
Y_x  = & - \, \displaystyle{\frac{\gamma_1 \, (\nu_2)_x}
{\gamma_2(\nu_1-\nu_2)}\, X}, & & \\
[4mm]
l_x  = & - \, \displaystyle{\frac{\nu_1 \, (\nu_2)_x}
{\gamma_2(\nu_1-\nu_2)}}\, X;  & &\\
[5mm]
X_y = & &\displaystyle{\frac{\gamma_2 \, (\nu_1)_y}
{\gamma_1(\nu_1-\nu_2)}\, Y,} &\\
[3mm]
Y_y  = &- \, \displaystyle{\frac{\gamma_2 \, (\nu_1)_y}
{\gamma_1(\nu_1-\nu_2)}}\, X & &
-\, \displaystyle{\frac{\nu_2 \,(\nu_1)_y}{\gamma_1
(\nu_1-\nu_2)}}\, l,\\
[4mm]
l_y  = &  & - \, \displaystyle{\frac{\nu_2 \, (\nu_1)_y}
{\gamma_1(\nu_1-\nu_2)}}\,Y. &
\end{array}$$
and the fundamental Bonnet theorem for strongly regular time-like surfaces states
as follows:
\begin{thm}\label{T:2.1}
Given four functions $\nu_1(x,y), \, \nu_2(x,y), \, \gamma_1(x,y),
\, \gamma_2(x,y)$ defined in a neighborhood $\mathcal D$ of $(x_0, y_0)$
satisfying the following conditions
$$\begin{array}{ll}
1) & \nu_1-\nu_2>0, \quad \gamma_1 \, (\nu_1)_y >0, \quad
\gamma_2 \, (\nu_2)_x > 0, \\
[4mm]
2.1) & \displaystyle{\left(\ln\frac{(\nu_1)_y}{\gamma_1}\right)_x=
\frac{(\nu_1)_x}{\nu_1-\nu_2},}
\qquad
\displaystyle{\left(\ln\frac{(\nu_2)_x}{\gamma_2}\right)_y=
-\frac{(\nu_2)_y}{\nu_1-\nu_2},}\\
[5mm]
2.2) & \displaystyle{\frac{\nu_1-\nu_2}{2}\left(\frac{(\gamma_2^2)_x}{(\nu_2)_x}-
\frac{(\gamma_1^2)_y}{(\nu_1)_y}\right)+ (\gamma_1^2+\gamma_2^2)=\nu_1\nu_2}
\end{array}$$
and an initial positive oriented orthonormal frame ${\bf z}_0X_0Y_0l_0$.

Then there exists a unique strongly regular time-like surface
${\mathcal M}: \; \textbf{z}=\textbf{z}(x, y), \; (x, y) \in \mathcal D_0 \;
((x_0, y_0) \in \mathcal D_0 \subset \mathcal D)$ with prescribed invariants
$\nu_1, \, \nu_2, \, \gamma_1, \, \gamma_2$ such that
$$z(x_0, y_0)=z_0, \; X(x_0, y_0)=X_0, \; Y(x_0, vy_0)=Y_0, \; l(ux_0, y_0)=l_0.$$
\end{thm}

Similarly to the case of minimal surfaces in Euclidean space we shall prove that any
minimal strongly regular surface admits locally geometric principal parameters. We can
assume that $\nu_1>0$ and we refer to the function $\nu=\nu_1$ as the normal
curvature function.
\begin{prop} Let $\mathcal M: \, {\bf z}={\bf z}(x,y), \; (x,y)\in \mathcal D$
be a minimal strongly regular surface whose parametric net is principal.
Then there exist locally principal parameters $(\bar x, \bar y)$ such that
$${\bf z}_{\bar x}^2={\bf z}_{\bar y}^2=\frac{1}{\nu}, \qquad \nu=\nu_1 > 0.$$
\end{prop}

{\it Proof:} Taking into account (2.1) and (2.3), we obtain
$$ (\ln \sqrt{\nu E})_y=0, \quad (\ln \sqrt{\nu G})_x=0,$$
which shows that $\nu E$ is only a function of $x$ and $\nu G$ is only a function of $y$.

Let $(x_0, y_0) \in \mathcal D$. We introduce new parameters $(\bar x, \bar y)$
in a neighborhood of $(x_0, y_0)$ by the formulas
$$\bar x = \int_{x_0}^x\sqrt{\nu E}\,dx +\bar x_0, \quad
\bar y = \int_{y_0}^y \sqrt{\nu G}\,dy +\bar y_0,$$
where $\bar x_0$ and $\bar y_0$ are constants. It follows immediately that
$(\bar x, \bar y)$ are again principal parameters and satisfy the required
property. $\qed$
\vskip 2mm
We call the parameters from the above proposition \emph{canonical principal} parameters.

Further we assume that the minimal strongly regular time-like surface $\mathcal M$
under consideration is parameterized with canonical principal parameters.

We use the following notations:
$$\nu =\nu_1>0, \quad \nu_2=-\nu<0, \quad \nu_1-\nu_2=2\nu>0.$$
Further we have
$$E=G=\frac{1}{\nu}, \quad e=-g=1, \quad \gamma_1=(\sqrt{\nu})_x, \quad \gamma_2=
-(\sqrt{\nu})_y.\leqno(2.5)$$

Then the equalities (2.2) become
$$\begin{array}{lccc}
X_x  = &  & \displaystyle{\frac{\nu_y}
{2\nu} \, Y} &- \sqrt {\nu} \, l,\\
[4mm]
Y_x  = & - \, \displaystyle{\frac{\nu_y}
{2\nu} \, X}, & & \\
[4mm]
l_x  = & - \, \sqrt{\nu} \, X;  & &\\
[5mm]
X_y = & &-\displaystyle{\frac{\nu_x}{2\nu} \, Y,} &\\
[3mm]
Y_y  = & \, \displaystyle{\frac{\nu_x}{2\nu}} \, X & &
+\,\sqrt{\nu} \, l,\\
[4mm]
l_y  = &  & + \, \sqrt{\nu} \,Y. &
\end{array}
\leqno(2.6)$$
and the integrability conditions of (2.6) reduce to
$$\Delta \ln \nu - 2 \nu=0. \leqno(2.7)$$

Theorem \ref{T:2.1} applied to minimal strongly regular time-like surfaces
parameterized with natural principal parameters implies:
\begin{cor}
Given a function $\nu (x,y) > 0$ in a neighborhood
$\mathcal D$ of $(x_0, y_0)$ with $\nu_x \nu_y \neq 0$,
satisfying the equation $(2.7)$ and an initial positive oriented orthonormal frame
${\bf z}_0X_0Y_0l_0$.

Then there exists a unique minimal strongly regular time-like surface
${\mathcal M}: \; \textbf{z}=\textbf{z}(x, y)$, \; $(x, y) \in \mathcal D_0 \;
((x_0, y_0) \in \mathcal D_0 \subset \mathcal D)$, for which

$(i)$ \; \; $(x, y)$ are canonical principal parameters;

$(ii)$ \, \, the invariants $\nu_1, \nu_2, \gamma_1, \gamma_2$ are the following
functions
$$\nu_1=\nu, \quad \nu_2= -\nu, \quad \gamma_1=(\sqrt {\nu})_y,  \quad
\gamma_2= -(\sqrt{\nu})_x;$$

$(iii)$ \, $z(x_0, y_0)=z_0, \; X(x_0, y_0)=X_0, \; Y(x_0, y_0)=Y_0, \; l(x_0, y_0)=l_0$.
\end{cor}
\vskip 2mm
Further we refer to (2.7) as the \emph{natural partial differential equation} of
minimal time-like surfaces.

The above statement gives a one-to-one correspondence between minimal strongly regular
time-like surfaces (considered up to a motion) and the solutions of the natural
partial differential equation, satisfying the conditions
$$\nu>0, \quad \nu_x \nu_y \neq 0. \leqno(2.8)$$

\section{Canonical Weierstrass representation of minimal strongly regular time-like surfaces}

Let $\mathcal M: \, \textbf{z}=\textbf{z}(x,y), \; (x,y) \in {\mathcal D}$ be a minimal
strongly regular time-like surface parameterized with canonical principal parameters.

Equalities (2.6) imply the following formulas for the Gauss map $l=l(x,y) \; (l^2=-1)$:
$$\begin{array}{ll}
l_{xx} &= \;\;\; \ds{\frac{\nu_x}{2\nu} \, l_x \, - \,
\frac{\nu_y}{2 \nu}\, l_y \, + \, \nu \, l,}\\
[4mm]
l_{xy} &= \;\;\; \ds{\frac{\nu_y}{2\nu}\, l_x \, + \, \frac{\nu_x}{2 \nu}\,l_y,}\\
[4mm]
l_{yy} &= - \ds{\frac{\nu_x}{2\nu}\,l_x \, + \, \frac{\nu_y}{2 \nu}\, l_y
\, + \,\nu \, l}
\end{array} \leqno(3.1)$$
and the normal vector function $l(x,y)$ satisfies the partial differential equation:
$$\Delta l-2\nu \, l=0.$$

The next statement makes precise the properties of the Gauss map of
a minimal time-like surface in terms of the canonical principal parameters.

\begin{prop}\label{P:3.1}
Let $\mathcal M: \textbf{z}=\textbf{z}(x,y), \; (x,y) \in {\mathcal D}$ be a minimal
strongly regular time-like surface parameterized with canonical principal parameters. Then the
Gauss map $l=l(x,y), \; (x,y)\in {\mathcal D}; \; l^2=-1$ has the following properties:
$$l_x^2=l_y^2=\nu >0, \quad l_x \,l_y=0,\quad \nu_x\nu_y\neq 0.\leqno(3.2)$$

Conversely, if a vector function $l(x,y), \; l^2=-1$ has the properties $(3.2)$, then there
exists locally a unique (up to a motion) minimal strongly regular time-like surface
$\mathcal M:\; \textbf{z}=\textbf{z}\,(x,y)$ determined by
$$\textbf{z}_x=-\frac{1}{\nu}\,l_x, \quad \textbf{z}_y=\frac{1}{\nu}\,l_y,
\leqno(3.3)$$
so that $(x,y)$ are geometrical principal parameters and $\nu(x,y)$ is the normal curvature
function of $\mathcal M$.
\end{prop}
\emph{Proof:} The equalities  $l_x=-\nu \, \textbf{z}_x, \; l_y=\nu \, \textbf{z}_y$,
and (2.5) imply (3.2).

For the inverse, it follows immediately that (3.2) implies (3.1). Furthermore, the second
equality of (3.1) implies that the system (3.3) is integrable.

Since (3.3) and (3.1) imply (2.6), it follows that
$\mathcal M:\;\textbf{z}=\textbf{z}(x,y)$ is a minimal strongly regular time-like surface
parameterized with canonical principal parameters with normal curvature function
$\nu=\,l_x^2=\,l_y^2 >0$.

Furthermore, it follows that the function $\nu=\,l_x^2=\,l_y^2$ satisfies the equation
(2.7). \qed
\vskip 2mm
Thus, any minimal strongly regular time-like surface locally is determined by the system (3.3),
where the vector function $l=l(x,y), \; l^2=-1$ satisfies the conditions (3.2). The so obtained
minimal surface is parameterized with canonical principal parameters.

\vskip 2mm
Now, let $H^2(-1):\xi^2 + \,\eta^2 - \,\zeta^2=-1$ be the unit
time-like sphere centered at the origin $O$ and \, $l(\xi, \eta,
\zeta), \zeta \neq 1$ be the position vector of an arbitrary point
on $H^2(-1)$, different from the pole $(0,0,1)$. Let us denote by $(x,y)$
the coordinates of any point in the parametric plane $\mathbb{C}$
and consider the standard conformal map $(x,y) \; \rightarrow \; (\xi, \eta, \zeta)$
of $H^2(-1)$ given by
$$l: \quad \begin{array}{l}
\displaystyle{\xi=\;\;\frac{2x}{x^2+y^2-1}},\\
[4mm]
\displaystyle{\eta=\;\;\frac{2y}{x^2+y^2-1}},\\
[4mm] \displaystyle{\zeta=-\frac{x^2+y^2+1}{x^2+y^2-1}};
\end{array}\qquad x^2+y^2-1\neq 0.$$
The vector function $l=l(x,y) \, (l^2=-1)$ has the properties
$$l_x^2=l_y^2=\frac{4}{(x^2+y^2-1)^2}, \quad l_x\,l_y=0$$
and generates the Enneper time-like surface.

The function
$$\nu=\frac{4}{(x^2+y^2-1)^2} \qquad x^2+y^2-1\neq 0$$
satisfies the equation (2.7).

Now we make more precise the Weierstrass representation of minimal strongly regular
time-like surfaces, stated in Theorem 1.
\vskip 4mm
\noindent
{\bf Proof of Theorem 1}
\vskip 4mm
Let $\mathcal M: \textbf{z}=\textbf{z}(x,y), (x,y) \in \mathcal D \subset \mathbb C$ be
a minimal strongly regular surface parameterized with canonical principal parameters.
Since the Gauss map $l=l(x,y) \; (l^2=-1$) of $\mathcal M$ is conformal and the orthogonal
frame field $l_x \, l_y \, l$ is left oriented, then the vector function $l$ is given
locally by
$$l: \quad \begin{array}{l}
\displaystyle{\xi=\;\;\frac{2u(x,y)}{u^2(x,y)+v^2(x,y)-1}},\\
[4mm]
\displaystyle{\eta=\;\;\frac{2v(x,y)}{u^2(x,y)+v^2(x,y)-1}},\\
[4mm]
\displaystyle{\zeta=-\frac{u^2(x,y)+v^2(x,y)+1}{u^2(x,y)+v^2(x,y)-1}};
\end{array}\qquad u^2(x,y)+v^2(x,y)-1\neq 0,\leqno(3.4)$$
where  $w=u(x,y)+i\,v(x,y)$ is a holomorphic function in $\mathbb C$.

We denote
$$\mu=\ds{\frac{(u_x^2+u_y^2)}{(u^2+v^2-1)^2}}.$$

If $\nu$ is the normal curvature function of $\mathcal M$, then the vector function
$\textbf{z}=\textbf{z}(x,y)$ satisfies the system
$$\begin{array}{l}
\textbf{z}_x=-\ds{\frac{1}{\nu}\,l_x=-\frac{1}{\nu}\,(u_x\,l_u+v_x\,l_v)},\\
[4mm]
\textbf{z}_y=\;\;\;\ds{\frac{1}{\nu}\,l_y=\;\;\;\frac{1}{\nu}\,(u_y\,l_u+v_y\,l_v)},
\end{array} \leqno(3.5)$$
which implies that
$\nu=\ds{\frac{4(u_x^2+u_y^2)}{(u^2+v^2-1)^2}}=4\mu$. Hence the holomorphic function $w$
satisfies the conditions
$$u^2+v^2-1 \neq 0; \qquad \mu>0, \qquad \mu_x \mu_y \neq 0.$$

Denoting by $w'=\ds{\frac{dw}{dz}=\frac{\partial w}{\partial z}=u_x-iu_y}$, we have
$$\frac{1}{w'}=\frac{u_x+iu_y}{u_x^2+v_y^2}$$
and (3.5) can be written in the form
$${\bf z}_x-i{\bf z}_y=-\frac{1}{w'} \frac{(u^2+v^2-1)^2}{4}(l_u+i l_v).\leqno(3.6)$$

Taking into account (3.4), we obtain from (3.6) the following formulas:
$$\begin{array}{l}
\ds{(z_1)_x-i(z_1)_y=\;\;\,\frac{1}{2} \, \frac{w^2+1}{w'}},\\
[4mm]
\ds{(z_2)_x-i(z_2)_y=-\frac{i}{2} \, \frac{w^2-1}{w'}},\\
[4mm]
\ds{(z_3)_x-i(z_3)_y=-\frac{w}{w'}},
\end{array}$$
which proves the assertion. $\hfill{\qed}$
\vskip 2mm
As an application we obtain a corollary for the solutions of the natural
partial differential equation.

\begin{cor} Any solution $\nu$ of the natural partial differential
equation $(2.7)$ satisfying the condition $(2.8)$ is given locally by the formula
$$\nu=\frac{4(u_x^2+u_y^2)}{(u^2+v^2-1)^2}, \quad  u^2+v^2-1\neq 0, \leqno(3.7)$$
where $w=u+iv$ is a holomorphic function in $\mathbb C$.
\end{cor}

{\it Proof:} Let $\nu(x,y)$ be a solution of (2.7) satisfying the conditions (2.8)
and let us consider the minimal strongly regular time-like surface $\mathcal M$, generated
by $\nu$. According to Theorem 1 it follows that the normal curvature function $\nu$
of $\mathcal M$ locally has the form (3.7).\qed

It is a direct verification that any function $\nu$, given by (3.7), where $w=u+iv \;
(u_x^2+u_y^2 > 0)$ is a holomorphic function in $\mathbb C$, satisfies (2.7).

\begin{rem}
The canonical Weierstrass representation of minimal strongly regular
time-like surfaces is based on the Gauss map of the Enneper time-like
surface $(w=z)$. Choosing the Gauss map of any other minimal strongly regular time-like
surface $\mathcal M$, we shall obtain its corresponding representation. This remark
is also valid for the representation $(3.7)$ of the solutions of the natural partial
differential equation $(2.7)$.
\end{rem}

\section{Maximal space-like surfaces and canonical principal parameters}

Let $\mathcal M: \, \textbf{z}=\textbf{z}(x,y), \; (x,y) \in {\mathcal D}$ be a space-like
surface in ${\R}_1^3$. In this case we assume that
$$E=\textbf{z}_x^2>0, \quad G=\textbf{z}_y^2<0, \quad  l^2=1.$$

We suppose that the surface has no umbilical points and the principal lines on
$\mathcal M$ form a parametric net, i.e.
$$F(x,y)=f(x,y)=0, \quad (x,y) \in \mathcal D.$$
Then the principal curvatures $\nu_1, \nu_2$ and the principal geodesic curvatures
(geodesic curvatures of the principal lines) $\gamma_1, \gamma_2$ are given by
$$\nu_1=\frac{e}{E}, \quad \nu_2=\frac{g}{G}; \qquad
\gamma_1=\frac{E_y}{2E\sqrt{-G}}, \quad \gamma_2= -\frac{G_x}{2\sqrt E \,G}. \leqno(4.1)$$

We consider the tangential frame field $\{X, Y\}$ determined by
$$X:=\frac{\textbf{z}_x}{\sqrt E}, \qquad Y:=\frac{\textbf{z}_y}{\sqrt{-G}}$$
and suppose that the moving frame $XYl$ is always positive oriented so that
$\nu_1-\nu_2 > 0$.

The following Frenet type formulas for the frame field $XYl$ are valid

$$\begin{tabular}{ll}
$\begin{array}{llccc}
\nabla_{X} \,X & = &  &\gamma_1 \,Y + \nu_1 \, l,  &\\
[2mm]
\nabla_{X} Y & = \;\;\;\gamma_1 \, X, & & \\
[2mm]
\nabla_{X} \, l & = - \nu_1 \, X, & & &
\end{array}$ &
\quad
$\begin{array}{llccc}
\nabla_{Y} \,X & = & & -\gamma_2 \, Y, &\\
[2mm]
\nabla_{Y} Y & = -\gamma_2 \, X & &  & - \nu_2 \, l,\\
[2mm]
\nabla_{Y}\, l & = & & - \nu_2 \, Y. &
\end{array}$
\end{tabular}\leqno(4.2)$$
\vskip 2mm

The Codazzi equations have the following form
$$\gamma_1=-\frac{Y(\nu_1)}{\nu_1-\nu_2}, \qquad
\gamma_2=-\frac{X(\nu_2)}{\nu_1-\nu_2} \leqno(4.3)$$
and the Gauss equation can be written as follows:
$$X(\gamma_2)+Y(\gamma_1)+ \gamma_1^2-\gamma_2^2 = \nu_1\nu_2 = K.)$$

A space-like surface ${\mathcal M}: \; \textbf{z}=\textbf{z}(x,y), \; (x,y)\in \mathcal D$
without umbilical points is said to be \emph{strongly regular} if (cf \cite{GM})
$$\gamma_1(x,y)\gamma_2(x,y) \neq 0, \quad (x,y)\in\mathcal D.$$

Since $$\gamma_1 \gamma_2 \neq 0 \; \iff \; (\nu_1)_y (\nu_2)_x \neq 0,$$
then the following formulas
$$\sqrt E=\frac{-(\nu_2)_x}{\gamma_2(\nu_1-\nu_2)} >0, \quad
\sqrt{-G}=\frac{-(\nu_1)_y}{\gamma_1 (\nu_1-\nu_2)}>0\, .\leqno(4.4)$$
are valid for strongly regular surfaces. Because of (4.4) formulas (4.2) become

$$\begin{array}{lccc}
X_x  = &  & \displaystyle{-\frac{\gamma_1 \, (\nu_2)_x}
{\gamma_2(\nu_1-\nu_2)}\, Y}
&- \, \displaystyle{\frac{\nu_1 \, (\nu_2)_x}{\gamma_2
(\nu_1-\nu_2)}\, l},\\
[4mm]
Y_x  = &\;\; \displaystyle{-\frac{\gamma_1 \, (\nu_2)_x}
{\gamma_2(\nu_1-\nu_2)}\, X}, & & \\
[4mm]
l_x  = & \, \;\;\;\displaystyle{\frac{\nu_1 \, (\nu_2)_x}
{\gamma_2(\nu_1-\nu_2)}}\, X;  & &\\
[5mm]
X_y = & &\displaystyle{\;\;\;\frac{\gamma_2 \, (\nu_1)_y}
{\gamma_1(\nu_1-\nu_2)}\, Y,} &\\
[3mm]
Y_y  = &\, \displaystyle{\;\;\;\frac{\gamma_2 \, (\nu_1)_y}
{\gamma_1(\nu_1-\nu_2)}}\, X & &
+\, \displaystyle{\frac{\nu_2 \,(\nu_1)_y}{\gamma_1
(\nu_1-\nu_2)}}\, l,\\
[4mm]
l_y  = &  & \;\;\;\displaystyle{\frac{\nu_2 \, (\nu_1)_y}
{\gamma_1(\nu_1-\nu_2)}}\,Y. &
\end{array}
$$
and the fundamental Bonnet theorem for strongly regular space-like surfaces states
as follows:
\begin{thm}\label{T:4.1}
Given four functions $\nu_1(x,y), \, \nu_2(x,y), \, \gamma_1(x,y),
\, \gamma_2(x,y); \; (x,y)\in \mathcal D$ satisfying the following
conditions:
$$\begin{array}{ll}
1) & \nu_1-\nu_2>0, \quad \gamma_1 \, (\nu_1)_y <0; \quad
\gamma_2 \, (\nu_2)_x < 0, \\
[4mm]
2.1) & \displaystyle{\left(\ln\frac{-(\nu_1)_y}{\gamma_1}\right)_x=
\frac{(\nu_1)_x}{\nu_1-\nu_2},}
\qquad
\displaystyle{\left(\ln\frac{-(\nu_2)_x}{\gamma_2}\right)_y=
-\frac{(\nu_2)_y}{\nu_1-\nu_2};}\\
[5mm]
2.2) & \displaystyle{\frac{\nu_1-\nu_2}{2}\left(\frac{(\gamma_1^2)_y}{(\nu_1)_y}+
\frac{(\gamma_2^2)_x}{(\nu_2)_x}\right)-\gamma_1^2+\gamma_2^2+\nu_1\nu_2=0.}
\end{array}$$
Then there exists a unique (up to a motion) strongly regular surface
with prescribed invariants
$\nu_1, \, \nu_2, \, \gamma_1, \, \gamma_2$.
\end{thm}
A space-like surface is maximal if $H=0$, i.e. $\nu_1+\nu_2=0$.
Next we introduce geometric principal parameters on maximal strongly regular space-like
surfaces. We can assume that $\nu_1>0$ and we refer to the function $\nu=\nu_1$ as the
normal curvature function.

\begin{prop} Let $\mathcal M: \, {\bf z}={\bf z}(x,y), \; (x,y)\in \mathcal D$
be a maximal strongly regular space-like surface whose parametric net is principal.
Then there exist locally principal parameters $(\bar x, \bar y)$ such that
$${\bf z}_{\bar x}^2=-{\bf z}_{\bar y}^2=\frac{1}{\nu}, \qquad \nu=\nu_1 > 0.$$
\end{prop}

{\it Proof:} Let $(x_0, y_0)$ be a point in $\mathcal D$. Taking into account (4.1)
and (4.3), we obtain
$$ (\ln \sqrt{\nu E})_y=0, \quad (\ln \sqrt{- \nu G})_x=0,$$
which shows that $\nu E$ is only a function of $x$ and $\nu G$ is only a function of $y$.
Introducing new parameters $(\bar x, \bar y)$ by the formulas
$$\bar x = \int_{x_0}^x\sqrt{\nu E}\,dx +\bar x_0, \quad
\bar y = \int_{y_0}^y \sqrt{- \nu G}\,dy +\bar y_0$$
we obtain that $(\bar x, \bar y)$ are again principal parameters and satisfy the required
property. $\qed$
\vskip 2mm
We call the parameters from the above lemma \emph{canonical principal} parameters.

From now on we assume that the maximal strongly regular space-like surface $\mathcal M$
is parameterized with canonical principal parameters.

We use the following notations:
$$\nu=\nu_1>0, \quad \nu_2=-\nu<0, \quad \nu_1-\nu_2=2\nu>0.$$
Further we have
$$E=-G=\frac{1}{\nu}, \quad e=g=1, \quad \gamma_1=-(\sqrt{\nu})_y, \quad \gamma_2=
(\sqrt{\nu})_x.\leqno(4.5)$$

Then the equalities (4.2) become
$$\begin{array}{lccc}
X_x  = &  & \displaystyle{-\frac{\nu_y}
{2\nu} \, Y} &+ \sqrt {\nu} \, l,\\
[4mm]
Y_x  = & - \, \displaystyle{\frac{\nu_y}
{2\nu} \, X}, & & \\
[4mm]
l_x  = & - \, \sqrt{\nu} \, X;  & &\\
[5mm]
X_y = & &-\displaystyle{\frac{\nu_x}{2\nu} \, Y,} &\\
[3mm]
Y_y  = & \, \displaystyle{-\frac{\nu_x}{2\nu}} \, X & &
+\,\sqrt{\nu} \, l,\\
[4mm]
l_y  = &  &  \, \sqrt{\nu} \,Y. &
\end{array}
\leqno(4.6)$$
and the integrability conditions of (4.6) reduce to
$$(\ln \nu)_{xx}-(\ln \nu)_{yy} + 2 \nu=0 \leqno(4.7)$$

Theorem \ref{T:4.1} applied to maximal strongly regular space-like surfaces
parameterized with canonical principal parameters implies:
\begin{cor}
Given a function $\nu (x,y) > 0$ in a neighborhood
$\mathcal D$ of $(x_0, y_0)$ with $\nu_x \nu_y \neq 0$,
satisfying the equation $(4.7)$ and an initial orthonormal frame ${\bf z}_0X_0Y_0l_0$.

Then there exists a unique strongly regular surface
${\mathcal M}: \; \textbf{z}=\textbf{z}(x, y), \; (x, y) \in \mathcal D_0 \;
((x_0, y_0) \in \mathcal D_0 \subset \mathcal D)$, such that

$(i)$ \; \; $(x, y)$ are canonical principal parameters;

$(ii)$ \, \, $z(x_0, y_0)=z_0, \; X(x_0, y_0)=X_0, \; Y(x_0, y_0)=Y_0, \; l(x_0, y_0)=l_0$;

$(iii)$ \, ${\mathcal M}$ is a maximal strongly regular space-like surface with invariants
$$\nu_1=\nu, \quad \nu_2= -\nu, \quad \gamma_1=-(\sqrt {\nu})_y,  \quad
\gamma_2=(\sqrt{\nu})_x.$$
\end{cor}
\vskip 2mm
Further we refer to (4.7) as the \emph{natural partial differential equation} of
minimal space-like surfaces.

The above statement gives a one-to-one correspondence between maximal strongly regular
space-like surfaces (considered up to a motion) and the solutions of the natural
partial differential equation, satisfying the condition

$$\nu>0, \quad \nu_x \nu_y \neq 0. \leqno(4.8)$$

\section{Canonical Weierstrass representation of maximal strongly regular space-like surfaces}

Let $H^2(1):\xi^2 + \,\eta^2 - \,\zeta^2=1$ be the unit
space-like sphere in ${\R}^3_1$ centered at the origin $O$ and \, $l(\xi, \eta,
\zeta)$ be the position vector of an arbitrary point on $H^2(1)$. If $\mathcal M: \,
\textbf{z}=\textbf{z}(x,y), \; (x,y) \in {\mathcal D}$ is a space-like surface, then its
Gauss map is $l: \; \mathcal D \longrightarrow H^2(1)$. Equalities (4.6) imply the following
formulas for the Gauss map:
$$\begin{array}{ll}
l_{xx} &= \ds{\frac{\nu_x}{2\nu} \, l_x \, + \,
\frac{\nu_y}{2 \nu}\, l_y \, - \, \nu \, l,}\\
[4mm]
l_{xy} &= \ds{\frac{\nu_y}{2\nu}\, l_x \, + \, \frac{\nu_x}{2 \nu}\,l_y,}\\
[4mm]
l_{yy} &= \ds{\frac{\nu_x}{2\nu}\,l_x \, + \, \frac{\nu_y}{2 \nu}\, l_y
\, +\,\nu \, l}
\end{array} \leqno(5.1)$$
The vector function $l(x,y), \; l^2=1$ satisfies the partial differential equation
$$l_{xx}-l_{yy}+2\nu \, l=0.$$

The next statement makes precise the properties of the Gauss map of
a maximal space-like surface in terms of canonical principal parameters.

\begin{prop}\label{P:3.1}
Let $\mathcal M: \textbf{z}=\textbf{z}(x,y), \; (x,y) \in {\mathcal D}$ be a maximal
strongly regular space-like surface parameterized with canonical principal parameters.
Then the Gauss map $l=l(x,y), \; l^2=1$ has the following properties:
$$l_x^2=-l_y^2=\nu >0, \quad l_x \,l_y=0,\quad \nu_x\nu_y\neq 0.\leqno(5.2)$$

Conversely, if a vector function $l(x,y), \, l^2=1$ has the properties $(5.2)$, then there
exists locally a unique (up to a motion) maximal strongly regular space-like surface
$\mathcal M:\; \textbf{z}=\textbf{z}\,(x,y)$ determined by
$$\textbf{z}_x=-\frac{1}{\nu}\,l_x, \quad \textbf{z}_y=\frac{1}{\nu}\,l_y,
\leqno(5.3)$$
so that $(x,y)$ are geometrical principal parameters and $\nu(x,y)$ is the normal
curvature function of $\mathcal M$.
\end{prop}
\emph{Proof:} The equalities  $l_x=-\nu \, \textbf{z}_x, \; l_y=\nu \, \textbf{z}_y$,
and (4.5) imply (5.2).

For the inverse, it follows immediately that (5.2) implies (5.1). Furthermore, the second
equality of (5.1) implies that the system (5.3) is integrable.

Since (5.3) and (5.1) imply (4.6), it follows that
$\mathcal M:\;\textbf{z}=\textbf{z}(x,y)$ is a maximal strongly regular space-like
surface parameterized with canonical principal parameters whose normal curvature function
is $\nu$.

Furthermore, it follows that the function $\nu=l_x^2=-l_y^2$ satisfies the equation
(4.7). \qed
\vskip 2mm
Thus, any maximal strongly regular space-like surface locally is determined by the system (5.3),
where the vector function $l=l(x,y), \; l^2=1$ satisfies the conditions (5.2). The so obtained
maximal surface is parameterized with canonical principal parameters.

\vskip 2mm
In the study of maximal space-like surfaces the role of the Gauss plane is played by the
Lorentz plane $\mathbb{L}=\{z=x+jy\}, \, j^2=1$. The metric in the Lorentz plane is given by
the quadratic form $x^2-y^2$.

Let us denote by $(x,y)$
the coordinates of any point in the parametric plane $\mathbb{L}$
and consider the parametrization $(x,y) \; \rightarrow \; (\xi, \eta, \zeta)$
of $H^2(1)$ given by
$$l: \quad \begin{array}{l}
\displaystyle{\xi=\frac{x^2-y^2-1}{x^2-y^2+1}},\\
[4mm]
\displaystyle{\eta=\frac{2x}{x^2-y^2+1}},\\
[4mm] \displaystyle{\zeta=\frac{2y}{x^2-y^2+1}},
\end{array}\qquad x^2-y^2+1\neq 0.\leqno(5.4)$$

This vector function $l(x,y), l^2=1$ has the properties
$$l_x^2=-l_y^2=\frac{4}{(x^2-y^2+1)^2}, \quad l_x\,l_y=0$$
and generates the space-like Enneper surface.

The map (5.4) is conformal and the function
$$\frac{4}{(x^2-y^2+1)^2} \qquad x^2-y^2+1\neq 0 $$
satisfies the equation (4.7).

Let $w=w(z)$ be a map in $\mathbb{L}$ given by
$$w=u+jv: \quad
\begin{array}{l}
u=u(x,y),\\
[2mm]
v=v(x,y).\end{array}
$$
The differential operators $\ds{\frac{\partial}{\partial z}}$ and
$\ds{\frac{\partial}{\partial \bar z}}$ are introduced as follows:
$$\frac{\partial}{\partial z}=\frac{1}{2}\left(\frac{\partial}{\partial x}+
j\,\frac{\partial}{\partial y}\right), \quad \frac{\partial}{\partial \bar z}=
\frac{1}{2}\left(\frac{\partial}{\partial x}-
j\,\frac{\partial}{\partial y}\right).$$
Then $w(z)$ is a holomorphic function in $\mathbb{L}$ if
$\ds{\frac{\partial w}{\partial \bar z}=0}$, i.e.
$$u_x=v_y, \quad u_y=v_x.$$
\vskip 2mm
Now we shall prove Theorem 2, which makes more precise the Weierstrass representation of maximal
strongly regular space-like surfaces.
\vskip 4mm
\noindent
{\bf Proof of Theorem 2}
\vskip 4mm
Let $\mathcal M: \textbf{z}=\textbf{z}(x,y), (x,y) \in \mathcal D \subset \mathbb{L}$ be
a maximal strongly regular space-like surface parameterized with canonical principal
parameters.
Since the Gauss map $l: \; \mathcal D \; \rightarrow \; H^2(1)$ of $\mathcal M$ is
conformal and the orthogonal frame field $l_x\,l_y\,l$ is negative oriented, then the
vector function $l$ is given locally by
$$l: \quad \begin{array}{l}
\displaystyle{\xi=\frac{u^2(x,y)-v^2(x,y)-1}{u^2(x,y)-v^2(x,y)+1}},\\
[5mm]
\displaystyle{\eta=\frac{2u(x,y)}{u^2(x,y)-v^2(x,y)+1}},\\
[5mm]
\displaystyle{\zeta=\frac{2v(x,y)}{u^2(x,y)-v^2(x,y)+1}},
\end{array} \qquad u^2(x,y)-v^2(x,y)+1\neq 0, \leqno(5.5)$$
where  $w=u(x,y)+j\,v(x,y)$ is a holomorphic function in $\mathbb L$.

We denote
$$\mu=\frac{u_x^2-u_y^2}{(u^2-v^2+1)^2}.$$
If $\nu$ is the normal curvature function of $\mathcal M$, then the vector function
$\textbf{z}=\textbf{z}(x,y)$ satisfies the system
$$\begin{array}{l}
\textbf{z}_x=-\ds{\frac{1}{\nu}\,l_x=-\frac{1}{\nu}\,(u_x\,l_u+v_x\,l_v)},\\
[4mm]
\textbf{z}_y=\;\;\;\ds{\frac{1}{\nu}\,l_y=\;\;\;\frac{1}{\nu}\,(u_y\,l_u+v_y\,l_v)},
\end{array} \leqno(5.6)$$
which implies that $\nu=\ds{\frac{4(u_x^2-u_y^2)}{(u^2-v^2+1)^2}}=4\mu$. Hence the holomorphic
function $w$ satisfies the conditions
$$u^2-v^2+1 \neq 0; \qquad \mu>0, \qquad \mu_x \mu_y \neq 0.$$

Denoting by
$w'=\ds{\frac{dw}{dz}=\frac{\partial w}{\partial z}=u_x+ju_y}$, we have
$$\frac{1}{w'}=\frac{u_x-ju_y}{u_x^2-v_y^2}.$$

Then (5.6) can be written in the form
$${\bf z}_x+j{\bf z}_y=-\frac{1}{w'}\frac{(u^2-v^2+1)^2}{4} (l_u-j l_v).$$

Taking into account (5.5), we obtain
$$\begin{array}{l}
\ds{(z_1)_x+j(z_1)_y=-\frac{w}{w'}}\\
[3mm]
\ds{(z_2)_x+j(z_2)_y=\frac{1}{2} \, \frac{w^2-1}{w'}},\\
[3mm]
\ds{(z_3)_x+j(z_3)_y=\frac{j}{2} \, \frac{w^2+1}{w'}},
\end{array}$$
which proves the assertion. $\hfill{\qed}$
\vskip 2mm
As an application we obtain a corollary for the solutions of the natural
partial differential equation.

\begin{cor} Any solution $\nu$ of the natural partial differential
equation $(4.7)$ satisfying the condition $(4.8)$ locally is given by the formula
$$\nu=\frac{4(u_x^2-u_y^2)}{(u^2-v^2+1)^2}, \quad  u^2-v^2+1 \neq 0, \leqno(5.8)$$
where $w=u+jv$ is a holomorphic function in $\mathbb L$.
\end{cor}

{\it Proof:} Let $\nu(x,y)$ be a solution of (4.7) satisfying the conditions (4.8)
and let us consider the minimal strongly regular space-like surface $\mathcal M$, generated
by $\nu$. According to Theorem 2 it follows that the normal curvature function $\nu$
of $\mathcal M$ locally has the form (5.8).\qed

It is a direct verification that any function $\nu$, given by (5.8), where $w=u+jv, \,
(u_x^2-u_y^2>0)$ is a holomorphic function in $\mathbb L$, satisfies (4.7).

\begin{rem} The canonical Weierstrass representation of maximal strongly regular space-like
surfaces is based on the Enneper space-like surface $(w=z)$. Choosing any other minimal
strongly regular space-like surface $\mathcal M$ as a basic surface, we shall obtain its
corresponding representation. This remark is also valid for the representation $(5.8)$ of
the solutions of the natural partial equation $(4.7)$.
\end{rem}

\end{document}